\documentclass[11pt,a4paper]{article}
\usepackage[utf8]{inputenc}
\usepackage{amsmath, amssymb, amsthm}
\usepackage{graphicx}
\usepackage{hyperref}
\usepackage{geometry}
\usepackage[percent]{overpic}
\usepackage{caption, subcaption}
\usepackage{amsthm, amssymb, amscd, amsmath, amsfonts, amsbsy}
\usepackage{stmaryrd}
\usepackage{float}
\usepackage{xfrac}
\usepackage{tikz-cd}
\usepackage{mathtools}
\usepackage{subcaption}
\usepackage{placeins}
\usepackage{longtable}
\usepackage{booktabs}
\usepackage{array}
\usepackage{xcolor}
\usepackage{authblk}

\newcommand{\reid}[1]{\raisebox{0.65em}{\;\;$\xleftrightarrow{\text{#1}}$\;\;}}
\newcommand\iconreid[1]{\raisebox{-0.7em}{\includegraphics[page=#1]{icons}}}

\title{Classification of (Uncolored) Bonded Knots and Links}

\author[1]{Matic Simonič}
\author[1,2,3]{Boštjan Gabrovšek}
\author[4]{Wanda Niemyska}

\affil[1]{University of Ljubljana, Faculty of Mathematics and Physics, Slovenia}
\affil[2]{Rudolfovo -- Science and Technology Centre Novo mesto, Slovenia}
\affil[3]{University of Ljubljana, Faculty of Education, Slovenia}
\affil[4]{University of Warsaw, Faculty of Mathematics, Informatics and Mechanics, Institute of Informatics, Poland}

\affil[ ]{\textit{Emails:} matic.simonic@fmf.uni-lj.si, bostjan.gabrovsek@pef.uni-lj.si, wanda@mimuw.edu.pl}

\date{\today}

\begin{document}

\maketitle

\begin{abstract}
We present a systematic classification of uncolored bonded knots with singularity number at most seven. 
Bonded knots provide a topological model for closed protein chains with intramolecular bridges, such as disulfide bonds. 
Following the tradition of knot tabulation, we describe a procedure based on the generation of planar graphs, their conversion into bonded knot diagrams, and the use of the Yamada polynomial together with brute-force Reidemeister moves to distinguish topological knotted types.
\end{abstract}

\section{Introduction}
\begin{figure}[t]
\centering
\includegraphics[width=0.9\linewidth]{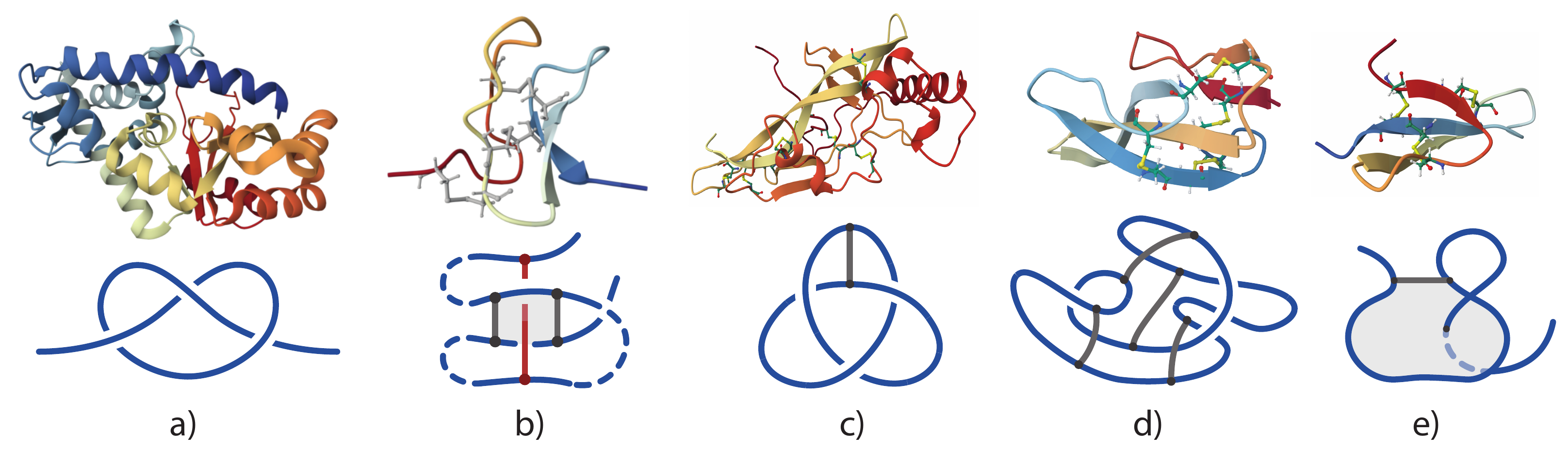}
\caption{
Representative examples of entangled protein topologies. 
Top: ribbon representations of protein structures; bottom: corresponding backbones with intramolecular bonds.
a) classical protein knot formed by the backbone chain, 
b) cystine motif where two bridges form a pierced disk, 
c) theta-curve formed by a single covalent bond connecting two segments of the backbone, 
d) bonded knot formed by multiple intramolecular bonds, and 
e) lasso protein topology, where a backbone segment pierces the minimal surface spanned by the backbone and a disulfide bond.
}
\label{fig:classes}
\end{figure}

The study of entanglement in biological polymers has become a vital area of research, as it is now known that approximately 1\% of entries in the Protein Data Bank (PDB) contain knots \cite{dabrowski2019knotprot} (see  \autoref{fig:classes}a). Although rare, knotted proteins occur in organisms across all domains of life, underscoring their fundamental role in biology \cite{perlinska2024everything}.
While classical knot theory focuses on closed loops, biological structures often include intramolecular bonds—such as disulfide bridges, hydrogen bonds, and salt bridges, which defines the molecule's folding pattern and stability \cite{sulkowska2008stabilizing}. Notably, when disulfide bridges are considered, roughly 30\% of proteins contain at least one such covalent bond \cite{RAJPAL20131721}, greatly expanding the percentage of biomolecules whose topological type can be analyzed through bonded entanglements.
These bonds give rise to a variety of topological motifs observed in proteins, including cystine knots, $\theta$-curve structures formed by a single bond, more complex bonded knots involving multiple bonds, and lasso proteins where a backbone segment pierces the surface spanned by the backbone and a bond (\autoref{fig:classes}).

Spatial graphs provide a natural model for such bonded conformations. 
In this setting, the polymer backbone is represented by a curve in space, 
while additional intramolecular bonds correspond to edges connecting 
different parts of the strand \cite{gabrovvsek2021invariant, gabrovsek2025bracket}. This perspective has proven useful in the 
study of DNA recombination and protein folding, where the presence of 
additional bonds significantly influences the global topology of the 
molecule.

To model these structures, we use the concept of \emph{bonded knots} \cite{gabrovvsek2021invariant, adams2020knot}.
This approach distinguishes between the protein backbone and the chemical bonds connecting different residues, extending the notion of ``circuit topology'' of proteins \cite{mashaghi2014circuit}, which captures only the combinatorial arrangement of bonds, to include the spatial knotting of the backbone.

In this paper, we classify bonded structures up to seven so called singularities (number of crossings plus number of vertices), extending previous work on $\theta$-curves (bonded knots with a single bond) performed by Moriuchi, who enumerated theta-curves with up to seven crossings in \cite{Moriuchi2009}.

\section{Definitions and Preliminaries}

A \emph{knotted graph} (or spatial graph) is a graph embedded in the 3-sphere $S^3$ \cite{Kauffman1989}.
A \emph{bonded knot} is defined as a pair $(L, \mathbf{b})$, where $L$ is a knot embedded in $S^3$, and $\mathbf{b} = \{b_1, b_2, \dots, b_n\}$ is a set of pairwise disjoint \emph{bonds} (closed intervals) properly embedded into $S^3 \setminus L$ such that their endpoints lie on $L$ and were introduced independently in \cite{adams2020knot}
 and \cite{gabrovvsek2021invariant}. In the uncolored case, we do not assign distinct chemical types to the bonds via a coloring function. Instead, we consider 3-valent graphs with a perfect matching, that is, there exists a subset of edges  incident to every vertex such that no two edges in the subset share a common vertex. Bonded knots and similar bonded structures have been studied in \cite{goundaroulis2017topological,gabrovvsek2021invariant,gugumcu2022invariants,gabrovsek2025bracket,diamantis2025topology,cavicchioli2025bonded}.

\subsection{Equivalence and Reidemeister Moves}
Two bonded knots are considered \emph{equivalent} if they are ambient isotopic, that is, they can be transformed into one another through a continuous deformation of the ambient space.

\begin{figure}[ht]
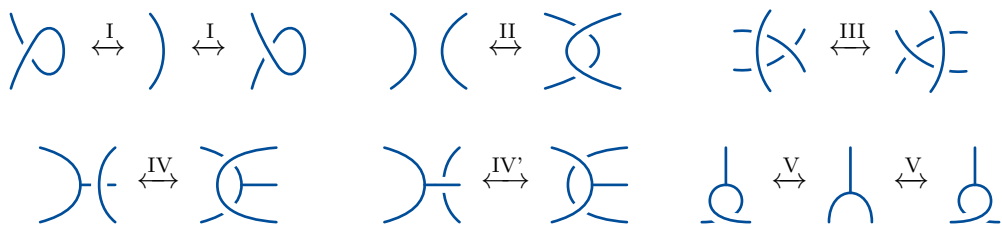

\centering
\begin{tabular}{ccccc}
$\iconreid{49}\reid{I}\iconreid{1}\reid{I}\iconreid{2}$ &\quad&
$\iconreid{3}\reid{II}\iconreid{4}$ &\quad&
$\iconreid{5}\reid{III}\iconreid{6}$ \\[2em]
$\iconreid{7}\reid{IV}\iconreid{8}$ &&
$\iconreid{9}\reid{IV'}\iconreid{54}$ &&
$\iconreid{55} \reid{V}\iconreid{56} \reid{V}\iconreid{57}$
\end{tabular}
% $$\iconreid{49}\reid{I}\iconreid{1}\reid{I}\iconreid{2} \hspace{4em}
%   \iconreid{3}\reid{II}\iconreid{4} \hspace{4em}
%   \iconreid{5}\reid{III}\iconreid{6} $$
% $$\iconreid{7}\reid{IV}\iconreid{8} \hspace{4em}
% \iconreid{9}\reid{IV'}\iconreid{54} \hspace{4em}
% \iconreid{55} \reid{V}\iconreid{56} \reid{V}\iconreid{57} $$ 
  \caption{Reidemesiter moves for 3-valent knotted graphs.} \label{fig:reidemeister}
  %The arcs in $I$, $II$, $III$, and the moving arc in $IV$ can also be (bolded) bonds.}
  \label{fig:reid}
\end{figure}

Diagrammatically, this equivalence is generated by a finite sequence of \emph{Reidemeister moves} presented in \autoref{fig:reid}. For bonded knots, these consist of the standard moves (I, II, and III), along with moves IV and V that handle interactions involving the vertices where bonds meet the backbone \cite{Kauffman1989, gabrovvsek2021invariant}. Kauffman \cite{Kauffman1989} distinguishes between rigid knotted graphs (without move V) and topological knotted graphs (with move V). Since rigid structures depend strongly on the chosen projection, we focus here only on the topological case, which is much harder to study due to lack of available invariants \cite{Kauffman1989}.

The \emph{crossing number} \(c(K)\) of a bonded knot \(K\) is the minimal number of crossings among all diagrams representing \(K\) under the generalized Reidemeister moves. 
Let \(b(K)\) denote the \emph{number of bonds} of \(K\). Since each bond corresponds to a pair of vertices, we define the \emph{singularity number}
\[
s(K) = c(K) + 2b(K).
\]
In this paper classify (uncolored) bonded knots with singularity number \(s(K) \le 7\).

\begin{figure} [ht]
    \centering
    \includegraphics[width=0.3\linewidth]{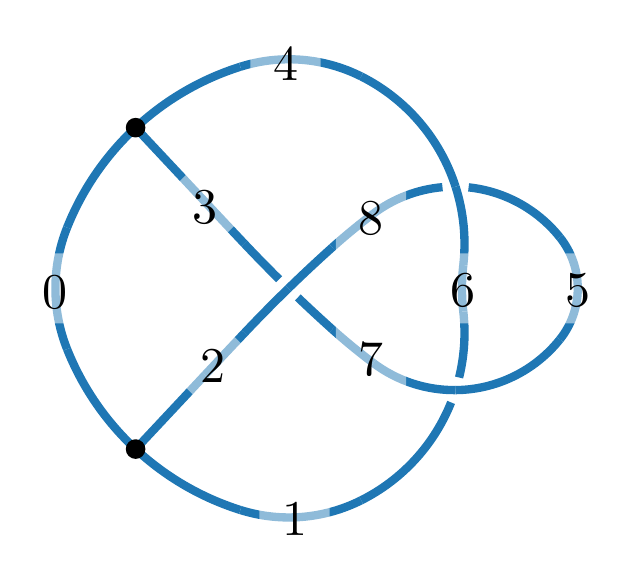}
    \caption{Enumeration of arcs of the bonded knot $B(3,1)_1.$}
    \label{fig:pd}
\end{figure}

A bonded knot can be encoded by \emph{Planar Diagram (PD)} notation as follows.
Let $K$ be a bonded knot from \autoref{fig:pd}. 
We enumerate the arcs consecutively from $0$ to $n-1$. 
For each vertex (bond endpoint or crossing), we record the incident arcs in counterclockwise order, beginning with an undercrossing arc (in case of a crossing).
The PD code is then the ordered list consisting of the two endpoint entries and the entries corresponding to all crossings. 
For the diagram shown in \autoref{fig:pd}, the PD code is
\[ \mathrm{PD}(B(3,1)_1) =\texttt{ 
V[0,1,2],V[0,3,4],X[1,4,5,6],X[7,8,3,2],X[8,9,10,5],X[6,10,9,7].}
\]

\section{Methodology}
The tabulation of knots dates back to the late nineteenth century. 
Motivated in part by Kelvin’s vortex atom theory, Kirkman and later Tait 
began compiling tables of knots organized by crossing number. Tait listed 
alternating knots with up to ten crossings, and Little later extended these 
tables to include non–alternating knots. Early tables were not entirely 
reliable; a famous example is the Perko pair, two diagrams originally listed 
as distinct ten–crossing knots that were later shown to represent the same 
knot. 

During the twentieth century, the development of knot invariants, beginning 
with the Alexander polynomial, provided systematic tools for distinguishing 
knots appearing in these tables. Later advances such as Conway notation and 
computer based methods using Dowker–Thistlethwaite codes enabled the 
enumeration of knots with much larger crossing numbers.

Another reason modern knot tabulation has become more tractable is the 
geometric structure of knot complements. By Thurston’s work, knot complements 
fall into three classes: torus, satellite, and hyperbolic. Most knots are hyperbolic, and the canonical triangulation of the hyperbolic complement can be computed efficiently. In practice, this provides one of the main tools for distinguishing knots with very large number of crossings.

In contrast, spatial graphs, including bonded knots, do not generally admit hyperbolic complements. Consequently, geometric tools used for classical knot tabulation are unavailable. Instead, classification relies on combinatorial invariants, particularly the Yamada polynomial, together with brute-force searches through generalized Reidemeister moves, making the tabulation problem significantly more difficult.

Our tabulation procedure consists of multiple steps, which we will describe in this section.
%is a multi-step computational process designed to systematically generate and distinguish bonded knot diagrams.
All computations were implemented in Python using the 
\texttt{KnotPy} \cite{knotpy} library, which provides tools for manipulating knot and 
spatial graph diagrams, performing Reidemeister moves, and computing 
invariants. The full source code used in the classification is available at \cite{simonic_gabrovsek_2025_bonded_knots}.

\begin{figure}[ht]
\centering

\begin{subfigure}{0.30\textwidth}
\centering
\includegraphics[page=1,width=\linewidth]{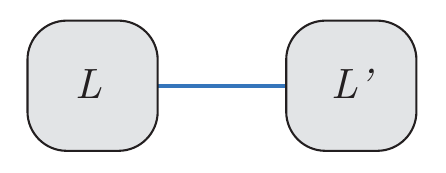}
\caption{$L \#_1 L'$} \label{fig:1}
\end{subfigure}
\hfill
\begin{subfigure}{0.30\textwidth}
\centering
\includegraphics[page=2,width=\linewidth]{sums.pdf}
\caption{$L \#_2 L'$}\label{fig:2}
\end{subfigure}
\hfill
\begin{subfigure}{0.30\textwidth}
\centering
\includegraphics[page=3,width=\linewidth]{sums.pdf}
\caption{$L \#_3 L'$}\label{fig:3}
\end{subfigure}

\caption{Connected sums of bonded knots. The operations $\#_1$, $\#_2$, and $\#_3$ correspond respectively to the order~1, order~2, and order~3 connected sums described by Moriuchi \cite{Moriuchi2009}.}
\end{figure}

\emph{Connected sums} of spatial graphs can occur in several different ways. 
Following Moriuchi \cite{Moriuchi2009}, an \emph{order~1 connected sum} corresponds to joining two diagrams by a bridge edge, producing a decomposition along a separating edge of the graph (\autoref{fig:1})
An \emph{order~2 connected sum} is analogous to the classical connected sum of knots, where two diagrams are joined along an arc of the backbone (\autoref{fig:2}).
Finally, Moriuchi introduced the notion of an \emph{order~3 connected sum}, in which two knotted graphs are joined at an 3-valent vertex (\autoref{fig:3}).

In this work we adopt the convention that a bonded knot is \emph{prime} if it cannot be decomposed as an order~1 or order~2 connected sum with two non-trivial knots (unknot) or trivial bonded knots (trivial theta curve or trivial handcuff link).
Order~3 connected sums are regarded as prime, since they do not admit a unique factorization.
Note that Moriuchi \cite{Moriuchi2009} regarded order 3 connected sums as non-prime.

The steps we used to classify prime bonded knots are the following.

\paragraph{Step 1: Generation.}
%This approach relies on the classical correspondence between knot diagrams 
%and $4$--regular planar graphs: crossings correspond to $4$–valent vertices 
%and edges represent the arcs of the diagram. By allowing $3$–valent vertices 
%that represent bond endpoints, planar graphs can be used to encode bonded 
%knot diagrams. 
We first generate all possible non-isomorphic graph diagrams up to seven vertices using the \texttt{plantri} software \cite{brinkmann2007fast}.
We choose connectivity $\geq 1$ graph, the minimum degree to $1$, the minimum number of edges to $\lceil n/2 \rceil$ (where $n$ denotes the number of vertices). The maximum number of edges is set to $2n$, corresponding to the number of edges of a $4$-regular graph on $n$ vertices.

We only consider graphs with vertices of degree $3$ and $4$. If a vertex has degree smaller than $4$, additional graphs are created by parallelizing the edges incident to that vertex until degree of a either vertex that the edge we would add is incident to would surpass $4$. In the resulting graphs, $4$-valent vertices are interpreted as crossings and $3$-valent vertices as endpoints of bonds. This interpretation imposes the additional requirement that the number of $3$-valent vertices must be even. Any graph not satisfying these conditions after the parallelization step is discarded. 

In total, this procedure produced $927$ candidate graphs.

\begin{figure}[ht]
    \centering
        \includegraphics[scale=1, page=13]{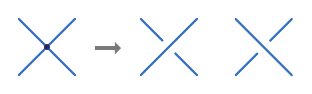}
    \caption{Replacing a 4-degree vertex with two possible crossing types.}
    \label{fig:crossing}
\end{figure}

\paragraph{Step 2: Conversion.}
Each $4$-valent vertex of the planar graphs is replaced by a crossing, under- and over-crossing, as in \autoref{fig:crossing} producing all possible bonded knot diagrams. This results in $20019$ candidate diagrams.

The diagrams are then simplified by applying brute-force Reidemeister moves and reducing them to canonical form using the KnotPy's \texttt{simplify} and \texttt{canonical} methods. The canonicalization converts a diagram into a unique  
representation that is invariant under diagrammatic symmetries.
The simplification procedure systematically performs all possible Reidemeister moves and explores sequences of such moves up to a fixed crossing-increasing depth and returns the minimal diagram. At this step we choose to reduce without crossing-increasing moves. In total, $1437$ diagrammatic canonical representatives remain. 

\paragraph{Step 3: Invariant Computation.}
For each diagram we compute the \emph{Yamada polynomial}, a polynomial 
invariant of spatial graphs introduced in \cite{Yamada1989}, and group the diagrams by this value.

Diagrams whose Yamada polynomial occurs only once form singleton groups and therefore represent diagrams that are distinct from all others by invariance of the polynomial. This produces $47$ such diagrams. The remaining $1390$ diagrams are partitioned into $137$ groups whose equivalence still needs to be determined.

\paragraph{Step 4: Simplification.}
For diagrams within the same group, we perform a brute-force comparison using sequences of Reidemeister moves implemented in KnotPy's \texttt{reduce\_equivalent\_diagrams} method. The algorithm explores the Reidemeister move spaces of all input diagrams simultaneously. Whenever the reachable sets of two diagrams intersect, the corresponding diagrams are identified as equivalent and only one (minimal) representative is retained; otherwise the diagram defines a distinct equivalence class.

The search is performed with increasing crossing-increasing depth in the Reidemeister move space. 

After depth $1$, we obtain $92$ additional unique diagrams, leaving $267$ undetermined diagrams distributed among $45$ groups.

After depth $2$, we obtain $13$ additional unique diagrams, leaving $128$ undetermined diagrams in $32$ groups.

After depth $3$, we obtain $2$ further unique diagrams, leaving $100$ undetermined diagrams in $30$ groups.

\begin{figure}[ht]
    \centering
        \includegraphics[scale=1, page=12]{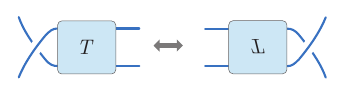}
    \caption{A flype move. While a flype can be decomposed into a sequence of Reidemeister moves, we treat it as an elementary move in some of our computations.}
    \label{fig:flype}
\end{figure}

Finally, depth $3$ with flypes included (\autoref{fig:flype}) produces no further reductions.
Since exploring depth $4$ becomes computationally expensive, the remaining $30$ groups are analyzed manually. Most groups consist of diagrams that are disjoint unions or connected sums, which can be identified at a glance. Only three groups remain that require further consideration, each containing two diagrams.

Two of these groups consist of pairs of mirror images, so only two groups need to be examined in detail. In the first group, the two diagrams are distinguished by removing the bond connecting the vertices: in one case this yields a Hopf link, while in the other it yields a link consisting of a trefoil knot and an unknot. In the second group, the diagrams differ in the number of components: one diagram has two components while the other has one.
\begin{figure}[htbp]
\centering

\begin{subfigure}{0.47\textwidth}
\centering
\includegraphics[width=0.44\linewidth,page=15]{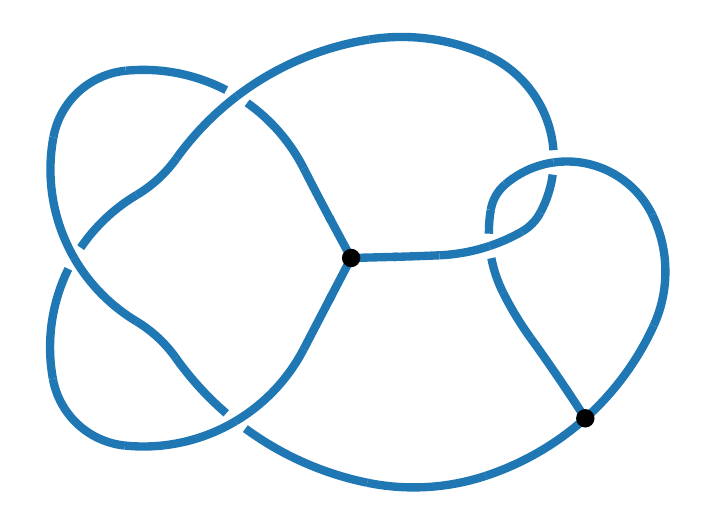}
\includegraphics[width=0.44\linewidth,page=16]{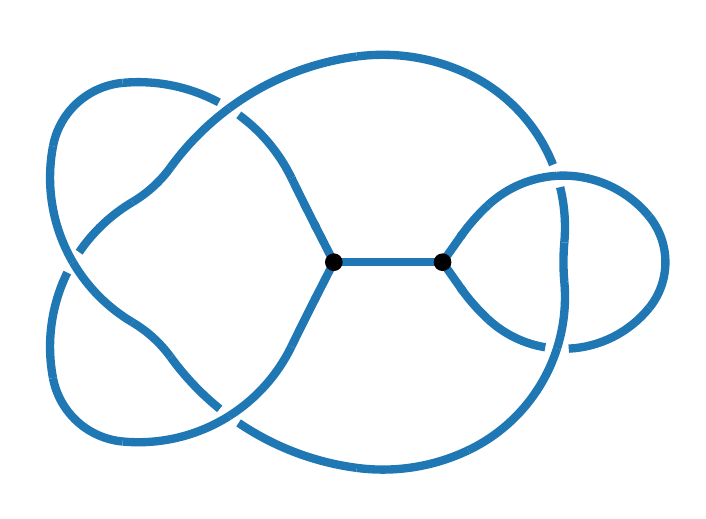}
\caption{First equivalence group}
\end{subfigure}
\hfill
\begin{subfigure}{0.47\textwidth}
\centering
\includegraphics[width=0.44\linewidth,page=8]{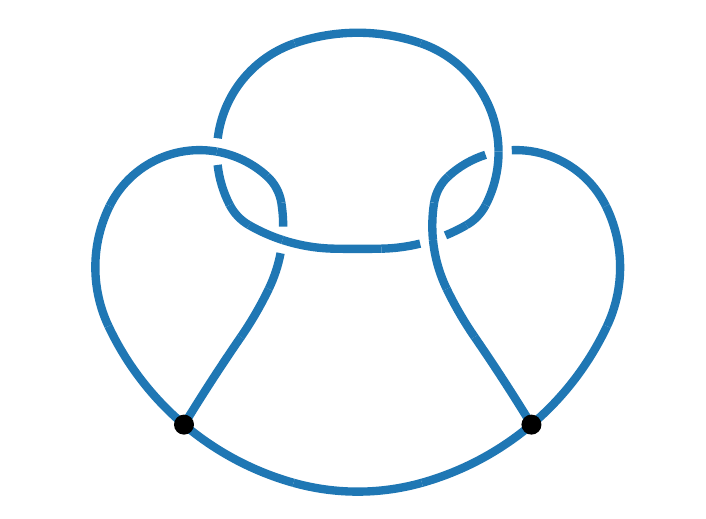}
\includegraphics[width=0.44\linewidth,page=7]{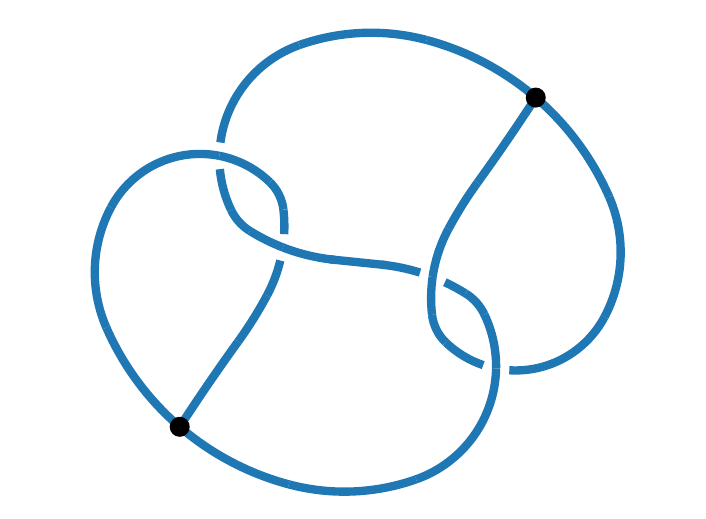}
\caption{Second equivalence group}
\end{subfigure}

\caption{The two remaining groups of (inequivalent) diagrams that have a shared Yamada polynomial and require manual analysis. They correspond to bonded links $H(5,1)_1$, $H(5,1)_2$, and $L(4,1)_1$, $H(4,1)_1$ in Tables \ref{fig:handcuff_knots} and \ref{fig:bonded_links}.}
\label{fig:remaining_groups}

\end{figure}
\paragraph{Step 5: Removal of non-admissible diagrams.}
At this stage we obtain $165$ unique graphs. We first remove diagrams corresponding to disjoint unions and connected sums, leaving $103$ diagrams. After removing classical knots and links, $51$ bonded knot diagrams remain.
Mirror images represent distinct bonded knot types, but the tables list only one representative from each mirror pair, with entries marked as chiral or achiral. To identify mirror pairs efficiently, diagrams are first partitioned according to their Yamada polynomials into those with unique and non-unique values.

For diagrams with unique Yamada polynomial, we compute the mirror image of each diagram and evaluate the Yamada polynomial of both the original and mirror diagrams. Since the list contains all possible diagrams, every polynomial from the original list appears in the mirror list. If the positions of the matching polynomials coincide, the diagram is achiral. Otherwise the diagram is chiral and the corresponding mirror diagram is removed. We are left with $31$ unique graphs, but as the diagram  in \autoref{fig:unbonded} is not a bonded knot we have $30$ final unique graphs.

%For diagrams with non-unique Yamada polynomials, the same procedure is applied, but additional care is required because several diagrams may share the same polynomial. In this case the correspondence between original and mirror diagrams is determined by their indices rather than polynomial values alone, after which the appropriate mirror duplicates are removed.

\section{Results}

Bonded knots obtained in our tabulation naturally split into three 
distinct classes according to the topological type of the underlying graph. 

\paragraph{Bonded knots.}
The class of bonded knots consists of those spatial graphs whose underlying graph does not 
contain a simple loop and is not a link (does not contain an extra connected component). These represent the generic case of bonded knots and we denote such diagrams by
\[
B(c,b)_i,
\]
where $c$ denotes the (minimal) number of crossings, $b$ denotes the number of bonds, 
and $i$ indicates the index among the bonded 
knots with the same values of $c$ and $b$. We have 20 such bonded knots.

\begin{figure}[ht]
\centering
\setlength{\tabcolsep}{8pt}
\renewcommand{\arraystretch}{1.2}
\begin{tabular}{cccc}
\includegraphics[page=2,width=0.22\textwidth]{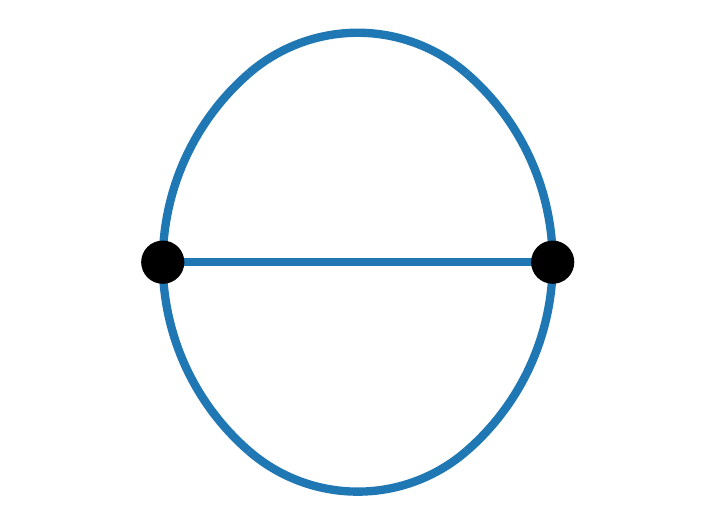} &
\includegraphics[page=5,width=0.22\textwidth]{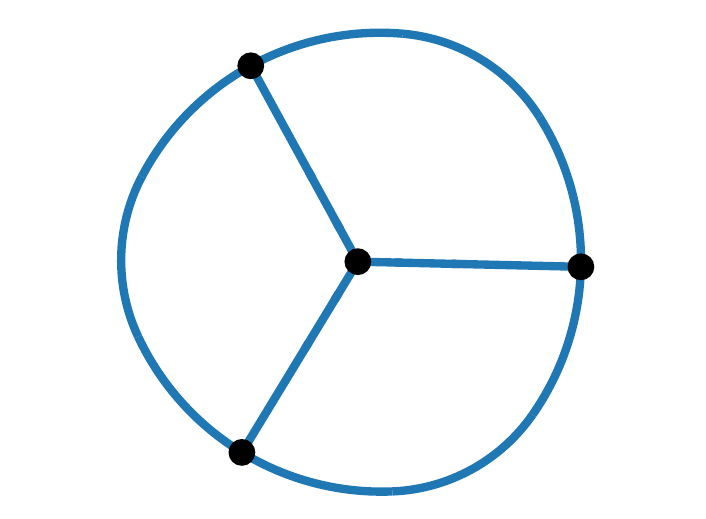} &
\includegraphics[page=6,width=0.22\textwidth]{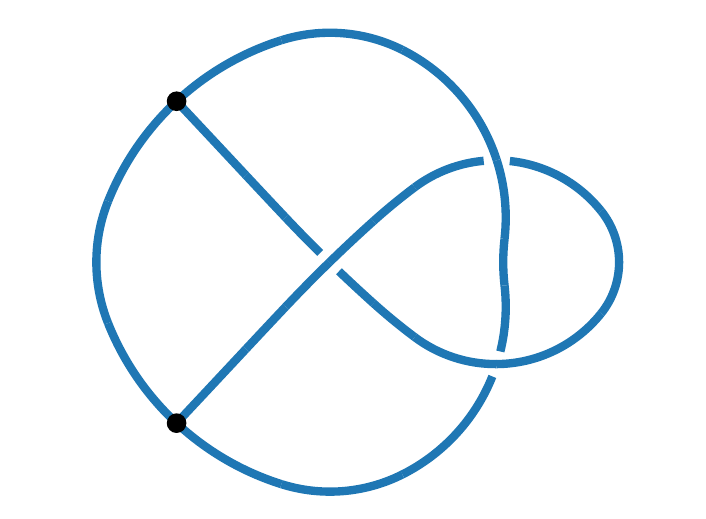} &
\includegraphics[page=9,width=0.22\textwidth]{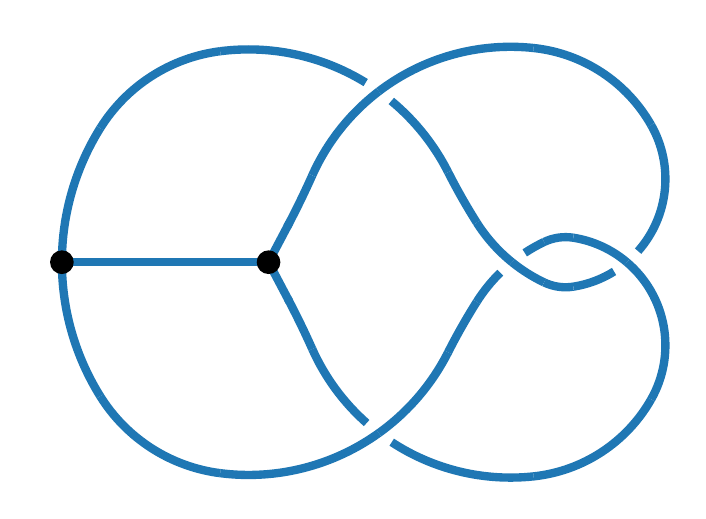} \\
$B(0,1)_1$ & $B(0,2)_1$ & $B(3,1)_1$ & $B(4,1)_1$ \\
\includegraphics[page=11,width=0.22\textwidth]{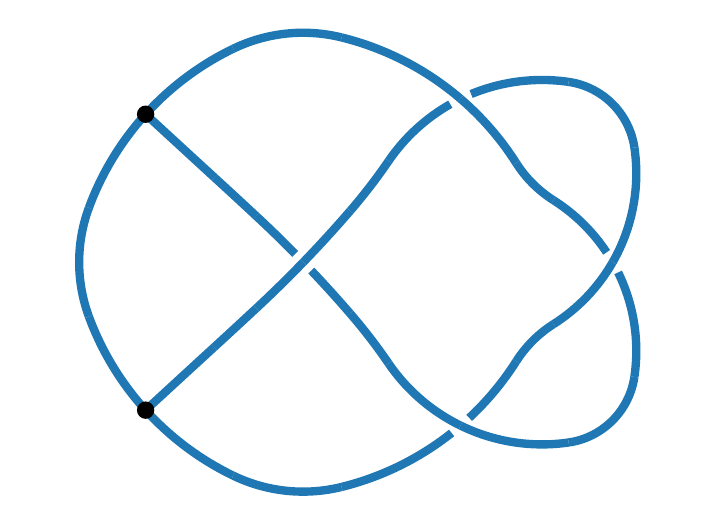} &
\includegraphics[page=13,width=0.22\textwidth]{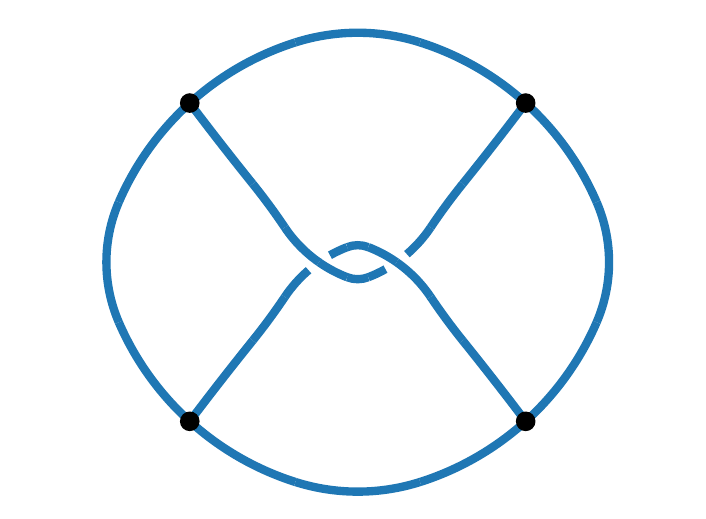} &
\includegraphics[page=14,width=0.22\textwidth]{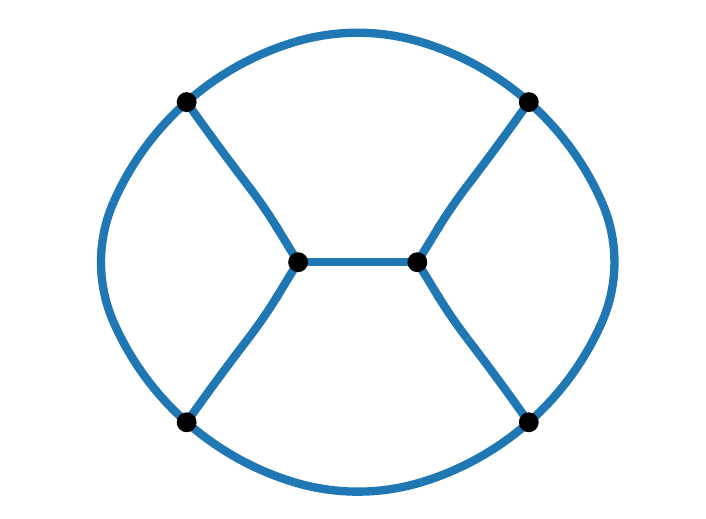} &
\includegraphics[page=17,width=0.22\textwidth]{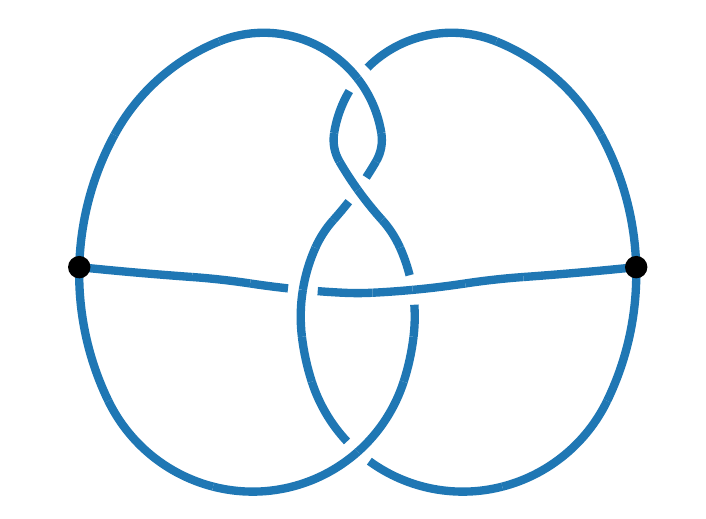} \\
$B(4,1)_2$ & $B(2,2)_1$ & $B(0,3)_1^{\#_3}$ & $B(5,1)_1$ \\
\includegraphics[page=18,width=0.22\textwidth]{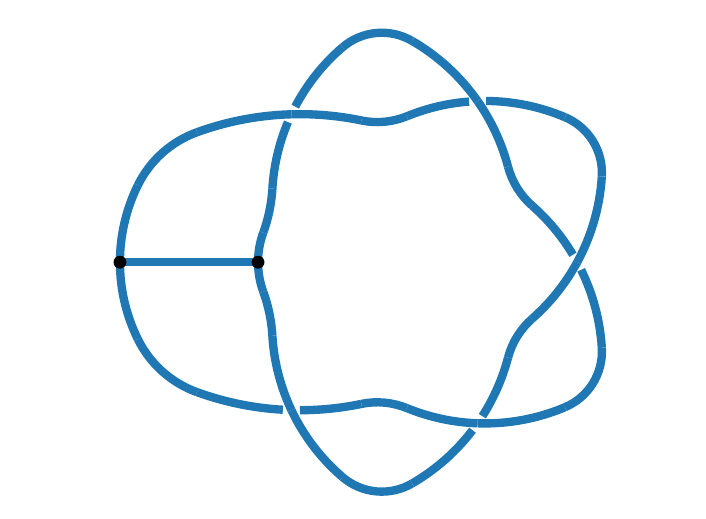} &
\includegraphics[page=20,width=0.22\textwidth]{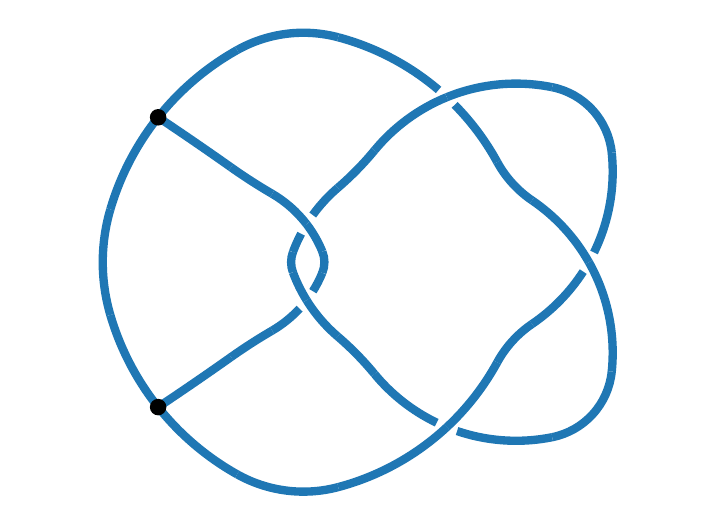} &
\includegraphics[page=21,width=0.22\textwidth]{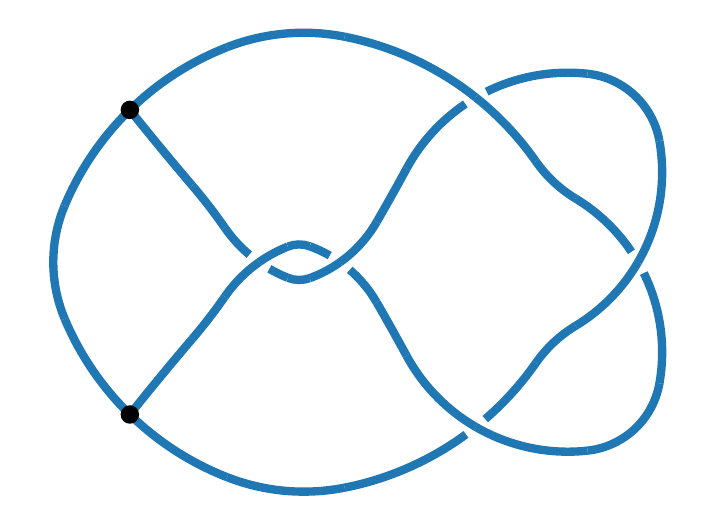} &
\includegraphics[page=24,width=0.22\textwidth]{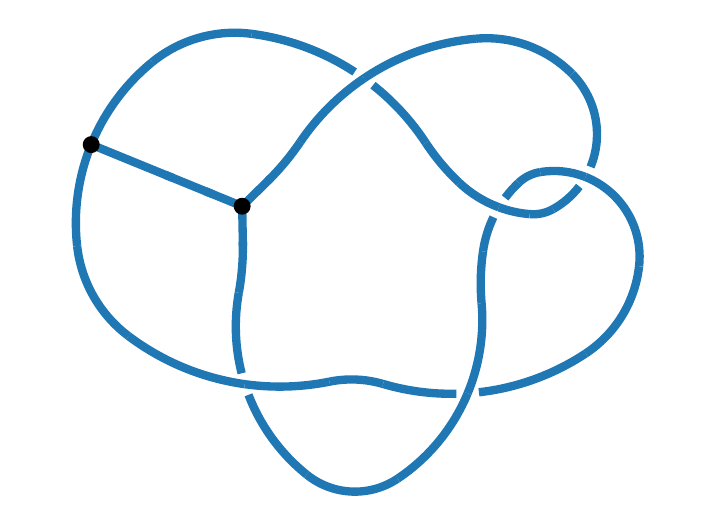}\\

$B(5,1)_2$ & $B(5,1)_3$ & $B(5,1)_4$ & $B(5,1)_5$  \\ 

\includegraphics[page=25,width=0.22\textwidth]{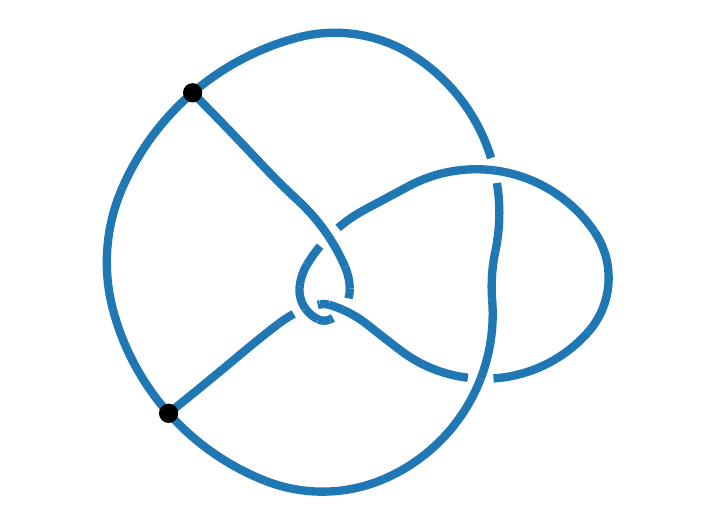}&
 \includegraphics[page=26,width=0.22\textwidth]{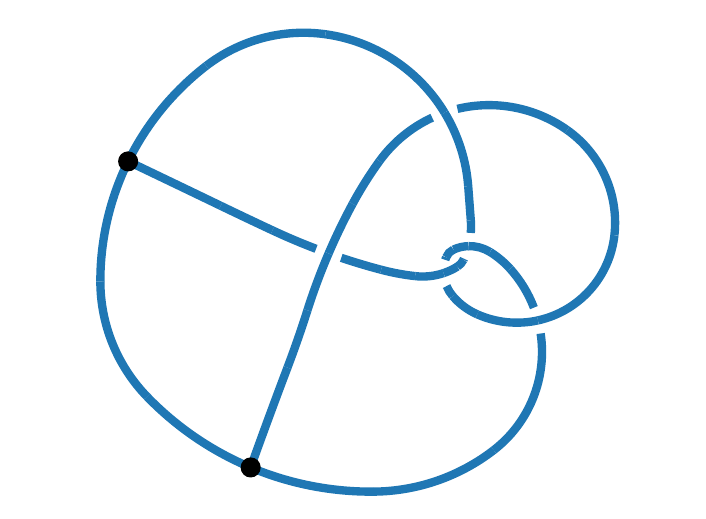}&
 \includegraphics[page=27,width=0.22\textwidth]{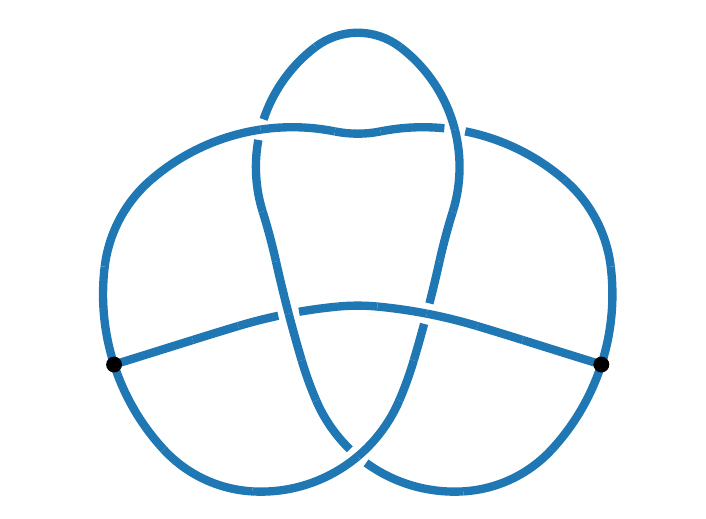}&
\includegraphics[page=22,width=0.22\textwidth]{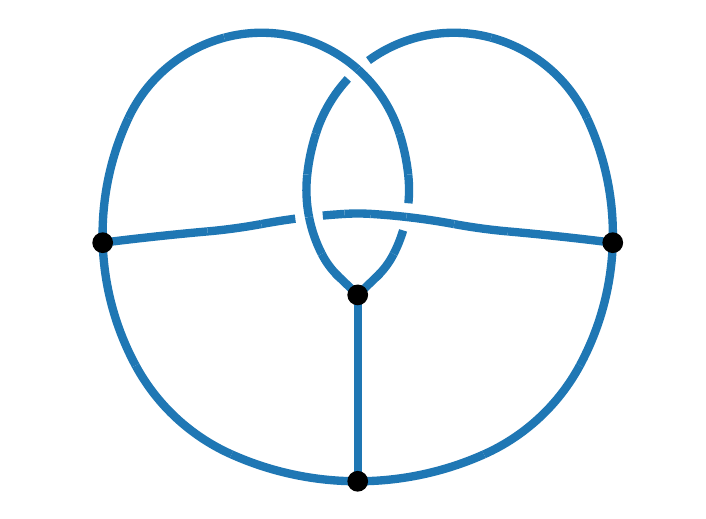}\\
   $B(5,1)_6$ & $B(5,1)_7$ & $B(5,1)_{8}$ & $B(3,2)_1$\\

\includegraphics[page=23,width=0.22\textwidth]{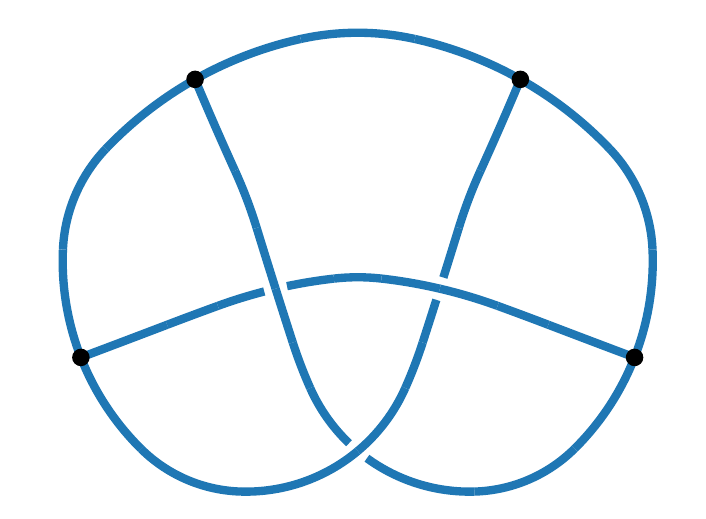}&
\includegraphics[page=28,width=0.22\textwidth]{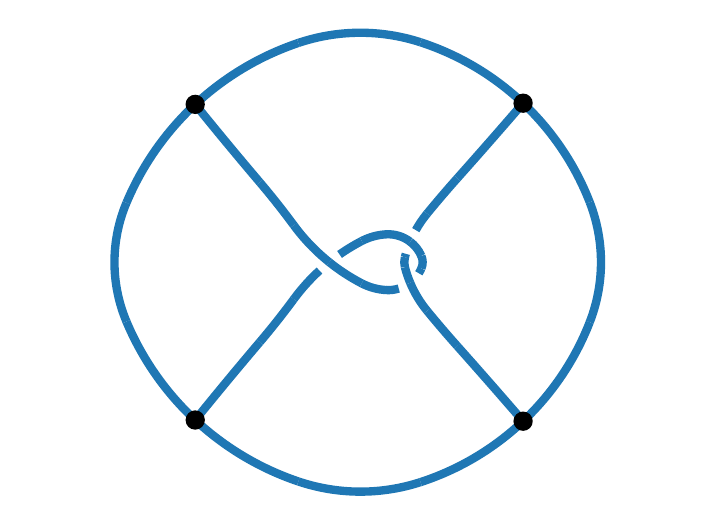} &
\includegraphics[page=30,width=0.22\textwidth]{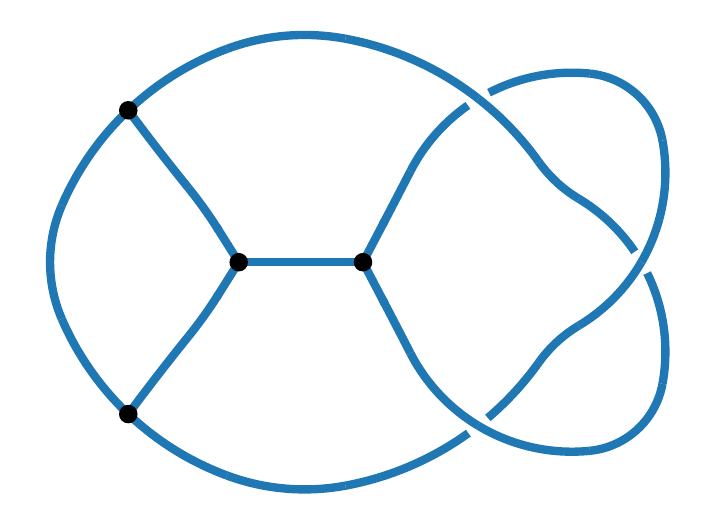} &
\includegraphics[page=31,width=0.22\textwidth]{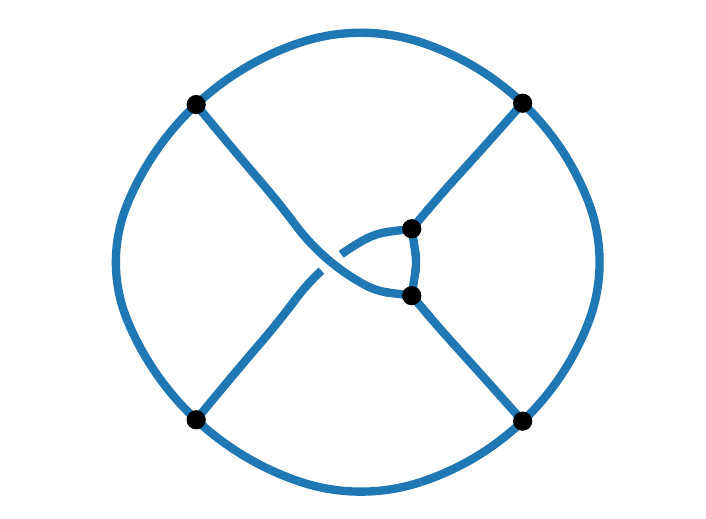} \\
 $B(3,2)_2$ & $B(3,2)_3$ & $B(3,2)_4^{\#_3}$ & $B(1,3)_1$ \\
\end{tabular}
\caption{Bonded knots. Each entry is represented by two rows: the diagram and its label.}
\label{fig:bonded_knots}
\end{figure}

In \autoref{fig:bonded_knots_yamada} the Yamada polynomials and chirality is presented, in \autoref{fig:bonded_knots_PD} the PD codes are presented, and in \autoref{fig:bonded_knots} the diagrams are depicted.

\paragraph{Bonded handcuff links.}
The class of bonded handcuff links consists of diagrams whose underlying graph 
contains a loop (a vertex adjacent to itself). These structures generalize \emph{handcuff links} and we denote these diagrams by
\[
H(c,b)_i,
\]
where again $c$ is the (minimal) number of crossings, $b$ is the number of bonds, and 
$i$ specifies the index. We have 6 such bonded knots.

\begin{figure}[ht]
\centering
\setlength{\tabcolsep}{8pt}
\renewcommand{\arraystretch}{1.2}
\begin{tabular}{cccc}
\includegraphics[page=1,width=0.22\textwidth]{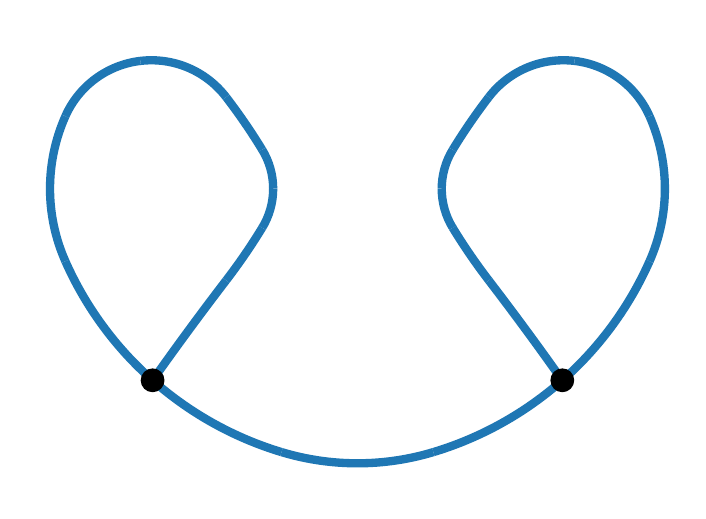} &
\includegraphics[page=4,width=0.22\textwidth]{bonded_knots.pdf} &
\includegraphics[page=7,width=0.22\textwidth]{bonded_knots.pdf} &
\includegraphics[page=12,width=0.22\textwidth]{bonded_knots.pdf} \\
$H(0,1)_1$ & $H(2,1)_1$ & $H(4,1)_1^{\#_3}$ & $H(2,2)_1^{\#_3}$ \\

\includegraphics[page=15,width=0.22\textwidth]{bonded_knots.pdf} &
\includegraphics[page=16,width=0.22\textwidth]{bonded_knots.pdf} &
 &  \\
$H(5,1)_1^{\#_3}$ & $H(5,1)_2^{\#_3}$ & & \\
\end{tabular}
\caption{Handcuff-like bonded knots. Each entry is represented by two rows: the diagram and its label.}
\label{fig:handcuff_knots}
\end{figure}

As in the previous case, \autoref{fig:handcuff_knots_yamada} lists the Yamada 
polynomials and chirality, \autoref{fig:handcuff_knots_PD} lists the PD codes, 
and \autoref{fig:handcuff_knots} depicts the corresponding diagrams.

\paragraph{Bonded links.}
The third class of bonded links consists of bonded diagrams with extra connected components. We denote these diagrams by
\[
L(c,b)_i,
\]
where again $c$ denotes the number of crossings, $b$ denotes the number of 
bonds, and $i$ iis the index. We have 4 such bonded links.

\begin{figure}[ht]
\centering
\setlength{\tabcolsep}{8pt}
\renewcommand{\arraystretch}{1.2}
\begin{tabular}{cccc}
\includegraphics[page=8,width=0.22\textwidth]{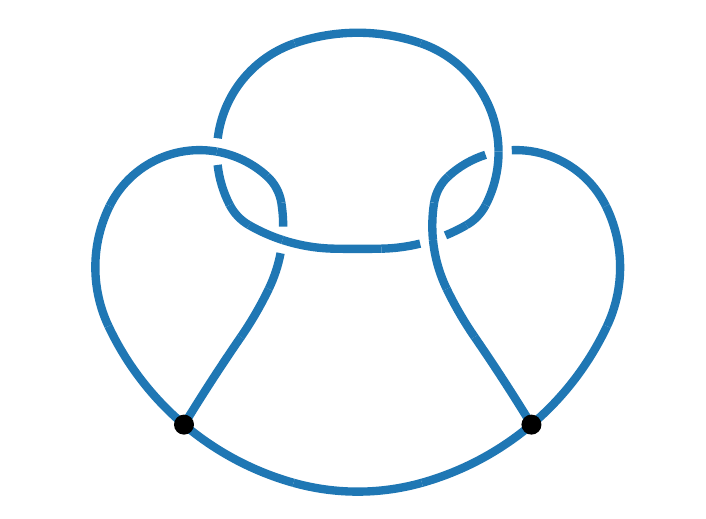} &
\includegraphics[page=10,width=0.22\textwidth]{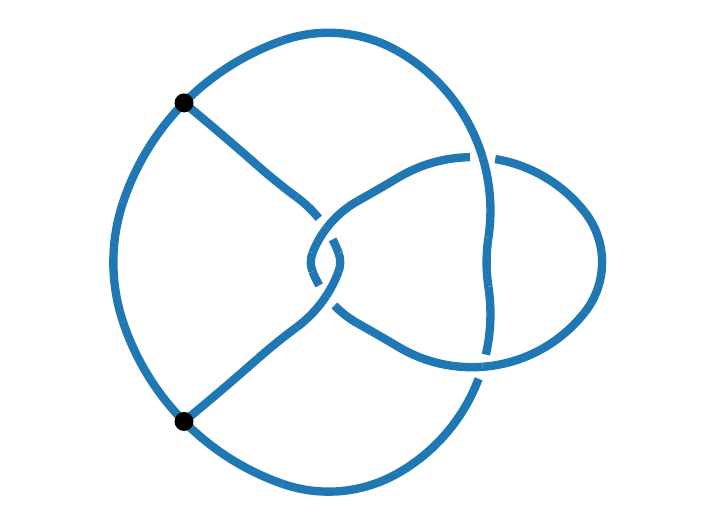} &
\includegraphics[page=19,width=0.22\textwidth]{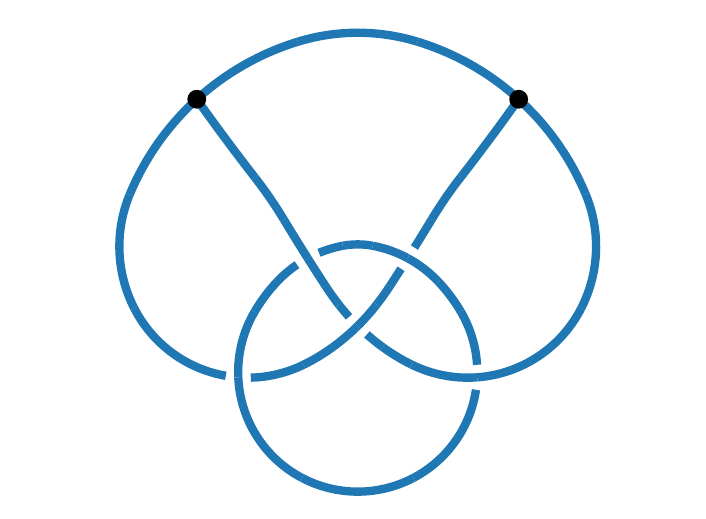} &
\includegraphics[page=29,width=0.22\textwidth]{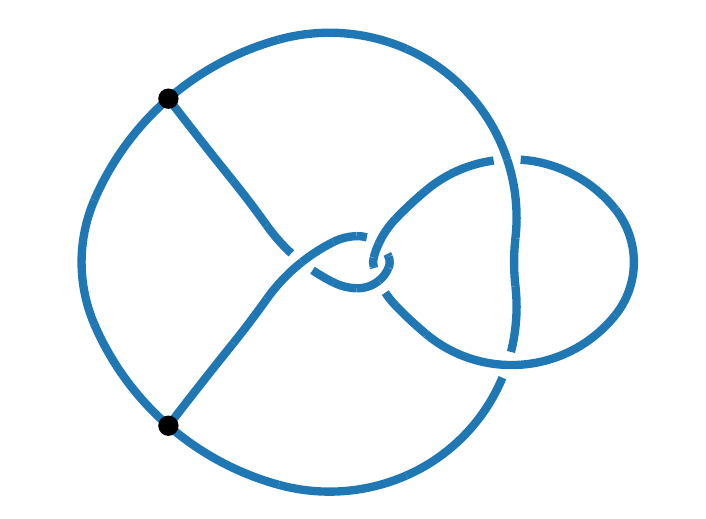} \\
$L(4,1)_1^{\#_3}$ & $L(4,1)_2$ & $L(5,1)_1$ & $L(5,1)_2$ \\
\end{tabular}
\caption{Linked bonded knots. Each entry is represented by two rows: the diagram and its label.}
\label{fig:bonded_links}
\end{figure}

As before, \autoref{fig:bonded_links_yamada} contains the Yamada polynomials and chirality, 
\autoref{fig:bonded_links_PD} lists the PD codes, and \autoref{fig:bonded_links} 
shows the corresponding diagrams.

%\paragraph{Connected sums.}
%Some diagrams appearing in the tables are obtained as connected sums of order 3, which we mark by a superscript $\#_3$. Across all three classes we identify 7 such bonded knots.

\paragraph{Summary of results.}
Up to seven singularities, there are 20 bonded knots, 6 bonded handcuff links, and 4 bonded links.  Across all three classes we identify 7 bonded knots that are connected sums of order 3 (marked by a superscript $\#_3$). There are 10 achiral bonded knots and 20 chiral bonded knots.
The Yamada polynomial is able to distinguish all but two cases of bonded knots (\autoref{fig:remaining_groups}), which are both order 3 connected sums of the same two bonded knots. Namely, the pair $(H(5,1)_1$ and $H(5,1)_2$ are both connected sums of order 3 of $B(3,1)_1$ and $H(2,1)_1$; the pair $L(4,1)_1$ and $H(4,1)_1$ are both connected sums of order 3 of two $H(2,1)_1$ handcuffs.

\begin{figure}[H]
\centering
\includegraphics[page=3,width=0.28\textwidth]{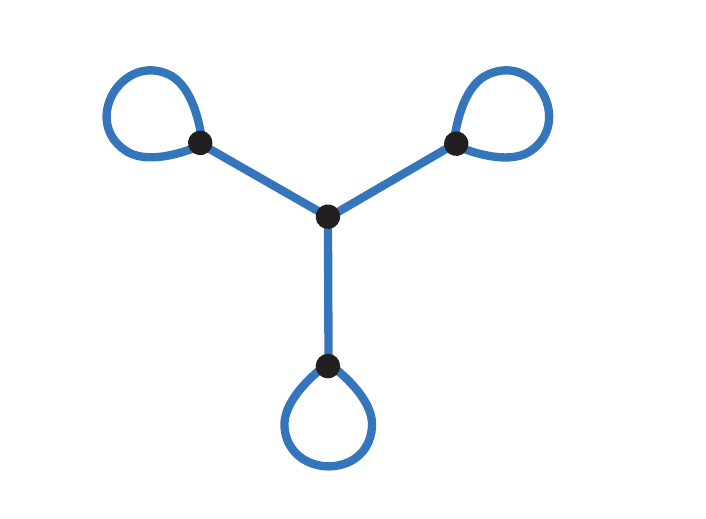}
\caption{The 3-valent graph, which is excluded from the tables.}\label{fig:unbonded}
\end{figure}

We obtain only one prime 3-valent knotted graph, which is not bonded, that is, the non-loop edges cannot be colored in such a way that two bonds do not meet. This is the (unknotted) graph depicted in \autoref{fig:unbonded}.

\begin{table}[h]
\centering
\begin{tabular}{lll}
\hline
Name & Chirality & Yamada \\
\hline
$B(0,1)_{1}$ & achiral & $-A^{4} - A^{3} - 2A^{2} - A - 1$ \\
$B(0,2)_{1}$ & achiral & $-A^{6} - 2A^{4} - 2A^{2} - 1$ \\
$B(3,1)_{1}$ & chiral & $-A^{12} - A^{11} - A^{10} - A^{9} - A^{8} - A^{6} - A^{4} + 1$ \\
$B(4,1)_{1}$ & chiral & $A^{15} + A^{12} + A^{9} + A^{7} + A^{5} + A^{4} + A^{2} - 1$ \\
$B(4,1)_{2}$ & chiral & $-A^{15} - A^{14} - A^{13} - A^{12} - A^{11} + A^{7} + A^{5} + A^{3} + A^{2} + 1$ \\
$B(2,2)_{1}$ & achiral & $A^{10} + A^{9} + A^{8} + 2A^{7} + 2A^{5} + 2A^{3} + A^{2} + A + 1$ \\
$B(0,3)_{1}^{\#_3}$ & achiral & $-A^{8} + A^{7} - 3A^{6} + 2A^{5} - 4A^{4} + 2A^{3} - 3A^{2} + A - 1$ \\
$B(5,1)_{1}$ & chiral & \resizebox{0.7\linewidth}{!}{$A^{17} - A^{16} - 2A^{15} + A^{14} - A^{13} + 2A^{11} + A^{10} + 2A^{9} + A^{7} - A^{5} + A^{4} + A^{3} - A^{2} + A + 1$} \\
$B(5,1)_{2}$ & chiral & $-A^{18} - A^{17} - A^{16} - A^{15} - A^{14} - A^{8} - A^{6} - A^{4} + A^{2} + 1$ \\
$B(5,1)_{3}$ & chiral & $-A^{17} + A^{16} + A^{13} + A^{11} + A^{9} - A^{7} - A^{6} - A^{5} - 2A^{4} - A^{3} - A^{2} - A - 1$ \\
$B(5,1)_{4}$ & chiral & $A^{18} + A^{15} - A^{13} - A^{11} - 2A^{10} - A^{8} + A^{7} - A^{6} - 2A^{3} - 1$ \\

$B(5,1)_{5}$ & chiral & $-A^{18} - A^{16} - 2A^{15} - A^{13} - 2A^{12} - A^{10} + A^{9} + A^{7} + A^{6} - A^{5} - A^{2} + 1$ \\
$B(5,1)_{6}$ & chiral & $2A^{17} + A^{16} + 2A^{14} + A^{13} + 2A^{11} + A^{9} - A^{8} - 2A^{5} - A^{2} + 1$ \\
$B(5,1)_{7}$ & chiral & $A^{18} - A^{16} + A^{15} - A^{13} + A^{12} - A^{10} - A^{8} - A^{6} + 2A^{4} - A^{3} + A^{2} + A - 1$ \\
$B(5,1)_{8}$ & chiral & $A^{16} + A^{15} - 2A^{11} - 2A^{10} - A^{9} - 2A^{8} - A^{6} + A^{5} + A^{2} - A - 1$ \\
$B(3,2)_{1}$ & chiral & $A^{13} - A^{9} - 3A^{8} - A^{7} - 4A^{6} - 2A^{4} + A^{3} + A^{2} + 2$ \\
$B(3,2)_{2}$ & chiral & $-A^{13} - A^{11} - A^{10} + A^{9} - A^{8} + 2A^{7} - A^{6} + A^{5} - A^{4} - A^{3} - A^{2} - A - 1$ \\
$B(3,2)_{3}$ & chiral & $-A^{13} + A^{12} + A^{9} - 2A^{8} - 3A^{6} - A^{5} - 2A^{4} - 2A^{3} - A^{2} - A - 1$ \\
$B(3,2)_{4}^{\#_3}$ & chiral & $-A^{14} - A^{12} - A^{11} - A^{10} - 2A^{8} + A^{7} - 2A^{6} + A^{5} - A^{4} + A^{2} - A + 1$ \\
$B(1,3)_{1}$ & achiral & $2A^{7} - A^{6} + 4A^{5} - 3A^{4} + 3A^{3} - 4A^{2} + A - 2$ \\
\hline
\end{tabular}
\caption{Bonded knots and their Yamada polynomials}
\label{fig:bonded_knots_yamada}
\end{table}

\begin{table}[h]
\centering
\begin{tabular}{lll}
\hline
Name & Chirality & Yamada \\
\hline
$H(0,1)_{1}$ & achiral & $0$ \\
$H(2,1)_{1}$ & achiral & $-A^{9} - A^{8} - A^{7} - A^{6} + A^{3} + A^{2} + A + 1$ \\
$H(4,1)_{1}^{\#_3}$ & achiral & $-A^{14} - A^{13} - A^{11} + A^{9} + A^{8} + 2A^{7} + A^{6} + A^{5} - A^{3} - A - 1$ \\
$H(2,2)_{1}^{\#_3}$ & achiral & $-A^{11} - A^{9} - A^{8} - A^{6} + A^{5} + A^{3} + A^{2} + 1$ \\
$H(5,1)_{1}^{\#_3}$ & chiral & $A^{17} - A^{15} + A^{14} - A^{13} - A^{12} - A^{10} - A^{8} - A^{5} + A^{4} + A^{3} + A + 1$ \\
$H(5,1)_{2}^{\#_3}$ & chiral & $A^{17} - A^{15} + A^{14} - A^{13} - A^{12} - A^{10} - A^{8} - A^{5} + A^{4} + A^{3} + A + 1$ \\
\hline
\end{tabular}
\caption{Handcuff-like bonded knots and their Yamada polynomials}
\label{fig:handcuff_knots_yamada}
\end{table}

\begin{table}[h]
\centering
\begin{tabular}{lll}
\hline
Name & Chirality & Yamada \\
\hline
$L(4,1)_{1}^{\#_3}$ & achiral & $-A^{14} - A^{13} - A^{11} + A^{9} + A^{8} + 2A^{7} + A^{6} + A^{5} - A^{3} - A - 1$ \\
$L(4,1)_{2}$ & chiral & $A^{14} + A^{13} + A^{12} + A^{11} + 2A^{10} + A^{9} + 2A^{8} + 2A^{7} + A^{6} + A^{5} + A^{3} + A^{2} + A + 2$ \\
$L(5,1)_{1}$ & chiral & \resizebox{0.8\linewidth}{!}{$-A^{17} + 2A^{15} + 3A^{12} + A^{11} + 2A^{10} + 3A^{9} + 2A^{8} + 2A^{7} + A^{6} + 2A^{5} + 2A^{4} - A^{3} + 2A^{2} - 2$} \\
$L(5,1)_{2}$ & chiral & \resizebox{0.8\linewidth}{!}{$A^{18} + 2A^{15} + A^{14} + 3A^{12} + 2A^{11} + 2A^{10} + 3A^{9} + 2A^{8} + 2A^{7} + A^{4} - 2A^{3} + A^{2} + A - 1$} \\
\hline
\end{tabular}
\caption{Linked bonded knots and their Yamada polynomials}
\label{fig:bonded_links_yamada}
\end{table}

\begin{table}[h]
\centering
\begin{tabular}{ll}
\hline
Name & PD code \\
\hline
$B(0,1)_{1}$ & {\scriptsize\ttfamily V[0,1,2],V[0,2,1]} \\
$B(0,2)_{1}$ & {\scriptsize\ttfamily V[0,1,2],V[0,3,4],V[1,4,5],V[2,5,3]} \\
$B(3,1)_{1}$ & {\scriptsize\ttfamily V[0,1,2],V[0,3,4],X[1,5,6,7],X[7,8,3,2],X[8,6,5,4]} \\
$B(4,1)_{1}$ & {\scriptsize\ttfamily V[0,1,2],V[0,3,4],X[1,4,5,6],X[7,8,3,2],X[8,9,10,5],X[6,10,9,7]} \\
$B(4,1)_{2}$ & {\scriptsize\ttfamily V[0,1,2],V[0,3,4],X[1,5,6,7],X[7,8,3,2],X[8,9,10,4],X[5,10,9,6]} \\
$B(2,2)_{1}$ & {\scriptsize\ttfamily V[0,1,2],V[0,3,4],V[1,5,6],X[2,7,8,3],V[4,9,5],X[9,8,7,6]} \\
$B(0,3)_{1}^{\#_3}$ & {\scriptsize\ttfamily V[0,1,2],V[0,3,4],V[1,5,6],V[2,7,3],V[4,8,5],V[6,8,7]} \\
$B(5,1)_{1}$ & {\scriptsize\ttfamily V[0,1,2],X[0,3,4,5],X[1,6,7,8],X[9,10,3,2],X[10,11,12,4],V[5,12,6],X[11,9,8,7]} \\
$B(5,1)_{2}$ & {\scriptsize\ttfamily V[0,1,2],X[0,3,4,5],V[1,5,6],X[6,7,8,2],X[3,9,10,4],X[7,11,12,8],X[9,12,11,10]} \\
$B(5,1)_{3}$ & {\scriptsize\ttfamily V[0,1,2],X[0,3,4,5],V[1,6,7],X[8,9,3,2],X[10,6,5,4],X[7,10,11,12],X[12,11,9,8]} \\
$B(5,1)_{4}$ & {\scriptsize\ttfamily V[0,1,2],V[0,3,4],X[1,5,6,7],X[8,9,3,2],X[10,11,12,4],X[5,12,11,6],X[7,10,9,8]} \\

$B(5,1)_{5}$ & {\scriptsize\ttfamily V[0,1,2],V[0,3,4],X[1,4,5,6],X[7,8,3,2],X[9,10,6,5],X[10,11,12,7],X[8,12,11,9]} \\
$B(5,1)_{6}$ & {\scriptsize\ttfamily V[0,1,2],X[0,3,4,5],V[1,6,7],X[8,9,3,2],X[10,6,5,4],X[7,11,12,8],X[9,12,11,10]} \\
$B(5,1)_{7}$ & {\scriptsize\ttfamily V[0,1,2],V[0,3,4],X[1,5,6,7],X[8,9,3,2],X[9,10,5,4],X[10,11,12,6],X[7,12,11,8]} \\
$B(5,1)_{8}$ & {\scriptsize\ttfamily V[0,1,2],X[0,3,4,5],X[1,6,7,8],X[8,9,3,2],X[9,10,11,4],V[5,11,12],X[6,12,10,7]} \\
$B(3,2)_{1}$ & {\scriptsize\ttfamily V[0,1,2],V[0,3,4],V[1,5,6],V[2,7,8],X[3,8,9,10],X[10,11,5,4],X[11,9,7,6]} \\
$B(3,2)_{2}$ & {\scriptsize\ttfamily V[0,1,2],V[0,3,4],V[1,5,6],X[2,7,8,3],X[4,8,9,10],V[5,10,11],X[11,9,7,6]} \\
$B(3,2)_{3}$ & {\scriptsize\ttfamily V[0,1,2],V[0,3,4],V[1,5,6],X[2,7,8,3],V[4,9,5],X[10,11,7,6],X[11,10,9,8]} \\
$B(3,2)_{4}^{\#_3}$ & {\scriptsize\ttfamily V[0,1,2],V[0,3,4],V[1,5,6],V[2,7,3],X[4,8,9,5],X[6,10,11,7],X[8,11,10,9]} \\
$B(1,3)_{1}$ & {\scriptsize\ttfamily V[0,1,2],V[0,3,4],V[1,5,6],X[2,7,8,3],V[4,9,5],V[6,10,7],V[8,10,9]} \\
\hline
\end{tabular}
\caption{PD codes of bonded knots}
\label{fig:bonded_knots_PD}
\end{table}

\begin{table}[h]
\centering
\begin{tabular}{ll}
\hline
Name & PD code \\
\hline
$H(0,1)_{1}$ & {\scriptsize\ttfamily V[0,0,1],V[1,2,2]} \\
$H(2,1)_{1}$ & {\scriptsize\ttfamily V[0,1,2],V[0,3,4],X[1,4,5,6],X[6,5,3,2]} \\
$H(4,1)_{1}^{\#_3}$ & {\scriptsize\ttfamily V[0,1,2],X[0,3,4,5],X[1,6,7,8],X[9,4,3,2],V[5,10,6],X[10,9,8,7]} \\
$H(2,2)_{1}^{\#_3}$ & {\scriptsize\ttfamily V[0,1,2],V[0,3,4],V[1,5,6],V[2,7,3],X[4,8,9,5],X[6,9,8,7]} \\
$H(5,1)_{1}^{\#_3}$ & {\scriptsize\ttfamily V[0,1,2],X[0,3,4,5],X[1,6,7,8],X[9,4,3,2],V[5,10,6],X[11,12,8,7],X[12,11,10,9]} \\
$H(5,1)_{2}^{\#_3}$ & {\scriptsize\ttfamily V[0,1,2],V[0,3,4],X[1,5,6,7],X[7,6,8,2],X[3,8,9,10],X[11,12,5,4],X[12,11,10,9]} \\
\hline
\end{tabular}
\caption{PD codes of handcuff-like bonded knots}
\label{fig:handcuff_knots_PD}
\end{table}

\begin{table}[h]
\centering
\begin{tabular}{ll}
\hline
Name & PD code \\
\hline
$L(4,1)_{1}^{\#_3}$ & {\scriptsize\ttfamily V[0,1,2],V[0,3,4],X[1,5,6,7],X[7,6,8,2],X[3,8,9,10],X[10,9,5,4]} \\
$L(4,1)_{2}$ & {\scriptsize\ttfamily V[0,1,2],V[0,3,4],X[1,5,6,7],X[7,8,9,2],X[3,9,8,10],X[10,6,5,4]} \\
$L(5,1)_{1}$ & {\scriptsize\ttfamily V[0,1,2],V[0,3,4],X[1,5,6,7],X[7,8,9,2],X[3,9,10,11],X[11,12,5,4],X[12,10,8,6]} \\
$L(5,1)_{2}$ & {\scriptsize\ttfamily V[0,1,2],V[0,3,4],X[1,5,6,7],X[8,9,3,2],X[10,6,5,4],X[7,11,12,8],X[9,12,11,10]} \\
\hline
\end{tabular}
\caption{PD codes of linked bonded knots}
\label{fig:bonded_links_PD}
\end{table}

%\input{bonded_tables}

%\input{bonded-knots-table}

% \begin{figure}
%     \centering
%     \includegraphics[width=0.75\linewidth]{bonded_sums.jpg}
%     \caption{Bonded links up to 7 vertices being connected sums of type $\#_3$}
%     \label{fig:placeholder}
% \end{figure}

\FloatBarrier
\subsection*{Acknowledgments}
B. Gabrovšek was financially supported by the Slovenian Research and Innovation Agency grants J1-4031, N1-0278, and program P1-0292. M. Simonič was financially supported by the Slovenian Research and Innovation Agency grant J1-4031 and program P1-0292. The authors declare that they have no conflict of interest.

%\section{Conclusion}
%We have established a classification of uncolored bonded knots up to seven crossings. This tabulation serves as a foundational resource for biological modeling, particularly in the study of protein circuit topology and the stabilizing role of intramolecular bridges [26, 27]. By extending the classification beyond simple $\theta$-curves, we provide the mathematical tools necessary to analyze the increasingly complex entangled structures being discovered in nature.

\bibliographystyle{abbrv}
\bibliography{biblio}
\end{document}